\documentclass[11pt]{article}
\usepackage{amssymb,amsmath}
\usepackage{graphicx,units,mathrsfs,bm}
\usepackage{subfigure}
\usepackage[sort&compress,numbers]{natbib} 
\usepackage{hyperref,xcolor}
\usepackage[paperwidth=8.5in, paperheight=11in, portrait, top=1in, bottom=1.in, left=1.15in, right=1.15in]{geometry}
\usepackage{tikz}
\usetikzlibrary{calc,decorations.markings,positioning}
\allowdisplaybreaks
\newcommand\bt{\raise 2pt \hbox{$\bigtriangledown$}\hskip 1.5pt}

\newmuskip\pFqskip
\mathchardef\pFcomma=\mathcode`, 

\numberwithin{equation}{section}

\newtheorem{prop}{Proposition}[section]

\newtheorem{defi}{Definition}[section]

\newtheorem{rem}{Remark}[section]

\newcommand{\BLUE}[1]{{\color{blue} #1}}

\begin{document}

\title{\bf 

An algebraic treatment of the Askey biorthogonal polynomials on the unit circle 
}
\author{

Luc Vinet\textsuperscript{$1$}\footnote{
E-mail: vinet@CRM.UMontreal.CA}~, 
Alexei Zhedanov\textsuperscript{$2$}\footnote{
E-mail: zhedanov@yahoo.com} \\[.5em]
\textsuperscript{$1$}\small~Centre de Recherches Math\'ematiques, 
Universit\'e de Montr\'eal, \\
\small~P.O. Box 6128, Centre-ville Station, Montr\'eal (Qu\'ebec), 
H3C 3J7, Canada.\\[.9em]

\textsuperscript{$2$}\small~School of Mathematics, 
Renmin University of China, Beijing, 100872, China
}
\date{\today}
\maketitle

\hrule
\begin{abstract}\noindent 
A joint algebraic interpretation of the biorthogonal Askey polynomials on the unit circle and of the orthogonal Jacobi polynomials is offered. It ties their bispectral properties to an algebra called the meta-Jacobi algebra $m\mathfrak{J}$.
\end{abstract}

\bigskip

\hrule

\section{Introduction}

In a commentary included in his edition of Szeg\H{o}'s collected works, Askey \cite{askey1982discussion} introduced sets of biorthogonal polynomials on the unit circle. These polynomials are defined as follows in terms of the standard Gauss hypergeometric series:
\begin{equation}
  P_n(z; \alpha, \beta) = \frac{(\beta)_n}{(\alpha+1)_n} \: {_2}F_1 \left( {-n, \alpha+1 \atop 1-\beta-n } ;z\right), \label{AP}   
\end{equation}
\begin{equation}
   Q_n(z; \alpha, \beta) =  P_n(z; \beta, \alpha), \label{part}
\end{equation}
where $(a)_k = a(a+1)\dots (a+k-1), k=1, 2, \dots \text{and} \;(a)_0=1$ are the Pochhammer symbols. 
The normalization is chosen so that $P_n(z)$ and $Q_n(z)$ are monic. The biorthogonality of these polynomials was proven in \cite{askey1985some} using slightly different conventions; it here reads:
\begin{equation}
   - \frac{1}{2\pi i}\oint_{|z|=1} dz \; {(-z)}^{-1-\beta} \left( 1-z \right) ^{\alpha+\beta} P_m(z, \alpha, \beta)\; Q_n(\frac{1}{z}, \alpha, \beta) = \frac{m!\;\Gamma(m+\alpha +\beta +1)}{\Gamma(\alpha +1) \Gamma(\beta +1)} \; \delta_{mn}, \label{bi}
\end{equation}
where $\Gamma(x)$ is the standard gamma function. Remember that $ (a)_n = \frac{\Gamma(n+ \alpha)}{\Gamma (\alpha)}$.
The branch of ${(-z)}^{-1-\beta}$ is chosen such that ${(-z)}^{-1-\beta}= |z|^{-1-\beta}$ if $arg \;z=\pi$ \cite{hendriksen1986orthogonal}, this is reflected in \cite{askey1982discussion}, \cite{askey1985some} by making the polar variable run from $-\pi$ to $\pi$.
 For the connection  with the spherical harmonics of the Heisenberg group, see \cite{greiner1983variations}, Sect 1.2. Let us also record that special cases of the Askey polynomials were obtained in \cite{shen2001orthogonal} as Fourier transforms of Laguerre polynomials (with weights attached). We refer to \cite{temme1986uniform} for historical remarks regarding these polynomials (see also \cite{borrego2016class}). 
 
 In his comments Askey expressed the opinion that the $P_n(z; \alpha, \beta)$ are the natural analogues of the Jacobi polynomials on the unit circle. We here reinforce this viewpoint by offering  a unified algebraic description of  these Askey polynomials on $S^1$ and of the Jacobi polynomials. This will involve the introduction of an algebra to be called \textit{meta-Jacobi} that will be seen to account for the bispectrality of both classes of functions.

Let us register for reference the definition and key properties of the monic Jacobi polynomials $\hat{P}_n^{(\alpha, \beta)}(x)$ defined on the interval $[0,1]$
\begin{equation}
   \hat{P}_n ^{(\alpha, \beta)}(x) = \frac{(-1)^n(1-\beta)_n} {(1 + \alpha + n)_n}\: {_2}F_1 \left( {-n, n + \alpha + 1 \atop 1 - \beta} ;x\right). \label{Jac}    
ƒ\end{equation}
Please note that for convenience an unconventional choice has been made for the parameters.
These polynomials possess the following orthogonality property
\begin{equation}
    \int_0^1  \hat{P}_m ^{(\alpha, \beta)}(x) \hat{P}_n ^{(\alpha, \beta)}(x)\; {x}^{-\beta} \left( 1-x \right) ^{\alpha+\beta}\; dx = h_n \;\delta_{mn} \qquad \alpha + \beta > -1, \quad \beta < 1, \label{orthrelJac}
\end{equation}
with the normalization factor $h_n$ given by
\begin{equation}
h_n= n! \; \frac{\Gamma(n - \beta +1) \Gamma(n + \alpha + 1) \Gamma(n + \alpha + \beta + 1)}{\Gamma (2n + \alpha + 1) \Gamma (2n + \alpha + 2)}.\label{hn}
\end{equation}
As is well known, in addition to satisfying a three-term recurrence relation, the polynomials $\hat{P}_n^{(\alpha, \beta)}(x)$ are eigenfunctions of the hypergeometric operator
\begin{equation}
    \mathcal{M} = x(x-1) \partial_x^2 + \left[ (\alpha+2)x + \beta-1   \right] \partial_x \label{M} 
\end{equation}
with eigenvalue $n(n + \alpha +1)$. These bispectral properties are encoded in the Jacobi algebra $\mathfrak{J}$ defined \cite{genest2016tridiagonalization} in terms of three generators $K_1$, $K_2$ and $K_3$ verifying the relations
\begin{align}
    [K_1, K_2] &= K_3 \label{K1} \\
    [K_2, K_3] &=aK_2^2 + bK_2, \label{K2}\\ 
    [K_3, K_1] &=a\{K_1,K_2\} + bK_1 + cK_2 +d, \label{K3}
\end{align}
where $[A, B]=AB-BA$, $\{A,B\}=AB+BA$ and $a, b, c, d$ are structure constants. Indeed $\mathfrak{J}$ is realized by taking 
\begin{equation}
    K_1=-\mathcal{M} \qquad K_2=x. 
\end{equation}
In this model where the generators $K_1$ and $K_2$ are the bispectral operators
\begin{equation}
    K_3= 2x(x-1) \partial_x^2 + \left[ (\alpha+2)x + \beta-1   \right] \partial_x, \label{K_3} 
\end{equation}
the parameters $a, b, c, d$ are
\begin{equation}
    a=2, \quad b=-2, \quad c=-\alpha (\alpha+2), \quad d=\alpha(1-\beta). \label{par}
\end{equation}

Headway in the algebraic description of bispectral biorthogonal functions was achieved recently by studying polynomial and rational functions of Hahn type \cite{tsujimoto2021algebraic}, \cite{vinet2020unified}. (Related Hahn rational functions also appear in \cite{koepf2007two} and \cite{masjed2012two}.) In broad strokes the general picture that emerges is as follows. Recall that generalized eigenvalue problems (GEVP) of the form $Md_n=\lambda_n Ld_n$ where $M$ and $L$ are two operators and $\lambda$ is the eigenvalue, naturally lead to biorthogonal functions which are rational (or polynomial) when $M$ and $L$ act tridiagonally in associated bases \cite{zhedanov1999biorthogonal}. Assume this to be the case. In the context mentioned above, it proved possible to adjoin a third operator $X$ to $M$ and $L$ such that the biorthogonal special functions are the overlaps between the relevant GEVP basis $\{d_n\}$ and the eigenbasis $\{e^*_z\}$ of the adjoint $X^T$ of $X$. As for the biorthogonal partner they are given reciprocally in terms of the bases for the corresponding adjoint problems. This offers a picture which is parallel to the description of hypergeometric (finite) polynomials using Leonard pairs \cite{terwilliger2001two}.
The differential/difference equation of the biorthogonal functions follows readily from the fact that $M-\lambda_n L$ which annihilates $d_n$ acts tridiagonally in the basis $\{e^*_z\}$. The second spectral equation stems from the observation that the operator $R^T=L^T X^T$ is such that $R^Te^*_{z} -z L^Te^*_{z} =0$
and that $R = XL$ acts tridiagonally on the basis $\{d_n\}$. The algebra generated by the triplet of operators $(M, L, R)$ which we have called the \textit{rational} Hahn algebra ($r\mathfrak{h}$) in the particular case treated in \cite{tsujimoto2021algebraic}, \cite{vinet2020unified} thus accounts for the two GEVPs that embody the bispectrality of the biorthogonal functions. Since $R$ factorizes as $XL$, the algebra generated by $(M, L, R)$ can be embedded in the \textit{meta}-algebra generated by $(M, L, X)$. The associated family of orthogonal polynomials also arises in this context as the overlaps between the eigenfunctions of the linear pencil $W=M + \mu L$ with the vectors $\{e^*_z\}$ (or equivalently as the scalar product of the eigenbases of the adjoint problems). The bispectrality of these polynomial functions is accounted for by the algebra generated by $(X, W)$. In our paradigm study, they are the Hahn polynomials, with $W$ and $X$ seen to generate the known Hahn algebra $\mathfrak{h}$.  In summary, for functions of the Hahn type, we observed that the meta-algebra $m\mathfrak{h}$ subsumes both $r\mathfrak{h}$ and $\mathfrak{h}$ and thus provides a unified description of both the biorthogonal and orthogonal families of functions. 
The two dimensional subalgebra of $m\mathfrak{h}$ generated by $M$ and $L$ is on its own remarkable since its three-diagonal representations lead alone to the corresponding orthogonal polynomials, that is the Hahn ones in this instance. The adjunction of $X$ to form the three-generated algebra has in fact the effect of constraining the representations of $M$ and $L$ to be three-diagonal in the eigenbasis of $X$.

We contend that this approach which allows the simultaneous description of hypergeometric orthogonal polynomials and associated families of biorthogonal functions extends beyond the Hahn functions case from which it is drawn. We shall add support to this suggestion by showing that the biorthogonal Askey polynomials on the unit circle together with the Jacobi polynomials are amenable to a unified treatment that follows the lines sketched above. In so doing we will provide an algebraic interpretation of the bispectral properties of the Askey polynomials  which is of interest in its own right. 
We might point out that it had been shown in \cite{grunbaum2004linear} that the recurrence relation of these polynomials can be obtained from a linear pencil in $\mathfrak{su}(1,1)$ without providing a full account of the bispectrality however.

The rest of the paper is organized as follows. The meta-Jacobi $m\mathfrak{J}$ is introduced and discussed in the next section. 
It is shown to be isomorphic to $\mathfrak{su}(1, 1)$.
The relevant generalized eigenvalue problems (GEVP) and eigenvalue problems (EVP) are solved on a $m\mathfrak{J}$ module and the appropriate overlaps are shown to yield the special functions of interest. The orthogonality relations are seen to follow from the completeness and orthogonality of the GEVP and EVP bases. The algebraic set-up is employed in Section \ref{4} to derive and interpret various properties of the Askey polynomials $ P_n(z; \alpha, \beta)$ and in particular their bispectrality. A differential model of $m\mathfrak{J}$ is obtained and used to obtain the differential equation and recurrence relation of  the polynomials $ P_n(z; \alpha, \beta)$ as well as some contiguity formulas. Perspectives are offered in the last section to conclude. Computational details are included in three appendices for completeness and the convenience of the reader.

\section{The meta-Jacobi algebra $m\mathfrak{J}$}\label{1}
The fundamental algebraic structure upon which the subsequent analysis hinges is introduced next.
\begin{defi}
The meta-Jacobi algebra $m\mathfrak{J}$ has generators $L$, $M$ and $X$ (and the central $1$)  verifying the commutation relations
\begin{align}
[L,M] &= L^2 -(\alpha+1)L -M \label{com_LM} \\
[L,X] &= X - 1 \label{com_LX} \\
[M,X] &= \{X,L\} -(\alpha+1) X + \beta \label{com_MX}. 
\end{align}
It is taken to be defined over the real numbers with the parameters $\alpha$ and $\beta$ in $\mathbb{R}$ unless specified otherwise.
\end{defi}
The Casimir element is checked to be
\begin{equation}
    Q= \{L^2,X\} -(\alpha + 1) \{L,X\} - \{M,X\} + 2 M + 2 \beta L.  \label{Casimir}
\end{equation}
We shall now observe that $m\mathfrak{J}$ is isomorphic to a Lie algebra. Recall that the Lie algebra $\mathfrak{su}(1,1)$ has the 
commutation relations:
\begin{equation}
    [J_0, J_{\pm}]=\pm J_{\pm}, \qquad [J_+, J_-]= -2 J_0, \label{su11}
\end{equation}
and the standard Casimir operator
\begin{equation}
    J^2 = J_0^2 - J_0 -J_+J_- . \label{Csu}
\end{equation}

We have:
\begin{prop} \label{iso}
The meta-Jacobi algebra $m\mathfrak{J}$ is isomorphic to the Lie algebra $\mathfrak{su}(1,1)$.
\end{prop}
This is confirmed by first observing that the commutation relations \eqref{su11} of $\mathfrak{su}(1,1)$ are recovered upon using the commutation relations \eqref{com_LM}, \eqref{com_LX}, \eqref{com_MX} of $m\mathfrak{J}$ and setting
\begin{align}
    J_0&=L - \frac{1}{2}(\alpha - \beta +1) \label{isom0} \\
    J_+&=X-1\label{isomplus} \\
    J_-&=-L^2 + (\alpha+1)L + M.\label{isomminus}
    \end{align}
    That we have an isomorphism is established by noting that this map is invertible and provides the following expressions of $L$, $M$ and $X$ in terms of the $\mathfrak{su}(1,1)$ generators:

\begin{align}
&L= J_0 +\frac{1}{2}(\alpha -\beta+1), \label{embed_slL} \\
&M= {J_0}^{2}+J_- - \beta \,J_0 -\frac{1}{4}(\alpha-\beta +1)(\alpha+\beta +1) \label{embed_slM} \\
&X = J_+ +1. \label{embed_slX} 
\end{align}
The isomorphism between the two-generated subalgebras spanned by $\{L, M\}$ and $\{J_0, J_-\}$ was observed in \cite{gaddis2015two}.
In light of the above formulas, the Casimir operator \eqref{Casimir} of the meta-Jacobi algebra can be expressed as
\begin{equation}
Q = 2 J^2 - \frac{1}{2} (\alpha - \beta +1)^2. \label{Q}
\end{equation}
\begin{rem}
In spite of this isomorphism it will be clear in the following that the $m\mathfrak{J}$ presentation is best suited for the algebraic interpretation of the Askey polynomials. We also stick to the terminology as it recalls the parallel with the treatment of the biorthogonal rational functions of Hahn type.
\end{rem}
\begin{prop} \label{embedding}
The Jacobi algebra $\mathfrak{J}$ defined in \eqref{K1}, \eqref{K2}, \eqref{K3} admits a simple embedding in the meta-Jacobi algebra $m\mathfrak{J}$. 
\end{prop}
This is seen by setting
\begin{equation}
    K_1=-M, \qquad K_2=X 
 \end{equation}
 and consequently
\begin{equation}
   K_3=-\{X,L\}+(\alpha+1)X-\beta. 
\end{equation}
Using the commutation relations \eqref{com_LM}, \eqref{com_LX}, \eqref{com_MX} of $m\mathfrak{J}$ it is straightforwardly verified that $K_1, K_2, K_3$ thus defined obey those of $\mathfrak{J}$ with the parameters given by
\begin{equation}
    a=2, \quad b=-2, \quad c=-\alpha(\alpha+2), \quad d= (\alpha+1)\beta -Q -1 . \label{d}
\end{equation}
Note the dependence of the parameter $d$ on the Casimir element $Q$. The distinctive feature of the meta-Jacobi algebra lies as we see in the fact that $K_3$ is resolved as a quadratic expression in terms of the fundamental generators $X$ and $L$.

\begin{rem}
In the following section we shall call upon representations of $\mathfrak{su}(1,1)$ and hence of $m\mathfrak{J}$ to interpret the Askey and Jacobi polynomials. In an irreducible representation, the Casimir element $J^2$ of $\mathfrak{su}(1,1)$ takes the form $\tau(\tau -1)$. Hereafter, we shall consider representations with
\begin{equation}
    \tau=\frac{1}{2}(\alpha + \beta +1).
\end{equation}
Equation \eqref{Q} which establishes the relation between the Casimir operator $Q$ of $m\mathfrak{J}$ and the one of $\mathfrak{su}(1,1)$ then yields for the value of $Q$:
\begin{equation}
    Q=2\alpha \beta -\alpha + \beta -1.
\end{equation}
Let us stress the coherence of the particular realization of the Jacobi algebra $\mathfrak{J}$ in terms of the bispectral operators of the Jacobi polynomials given in the Introduction with the embedding of $\mathfrak{J}$ in $m\mathfrak{J}$ given in Proposition \ref{embedding}. Indeed we see that with these choices for the Casimir elements, the parameter $d$ of the Jacobi algebra as given in \eqref{d} takes the proper value: $d=(\alpha +1)\beta-Q-1=\alpha(1-\beta)$.
\end{rem}


\section{Representations of the meta-Jacobi algebra and special functions}\label{3}
In this section we shall establish the connection between the meta-Jacobi algebra $m\mathfrak{J}$, the Askey polynomials  $P_n(z; \alpha, \beta)$, their biorthogonal partners $Q_n(z; \alpha, \beta)$ and the Jacobi polynomials $\hat{P}_n^{(\alpha, \beta)}(x)$. 
To that end we shall consider a $m\mathfrak{J}$ representation space inferred from the isomorphism of this algebra with $\mathfrak{su}(1,1)$.
We shall obtain the bases associated to the various EVP and GEVP defined on the chosen module to show that their overlaps are essentially the special functions mentioned above. This will cast these functions in their proper algebraic framework and readily lead to their (bi)orthogonality relations. We shall be working on a real infinite dimensional space equipped with a scalar product denoted by $\langle\;|\;\rangle$. $A^T$ will stand for the transpose of $A$ : $(\langle u|A^T) |v\rangle = \langle u| ( A|v\rangle)$.

Consider the infinite-dimensional module $\mathfrak{V}(\tau)$
with $\tau 
\in \mathbb{R}$ 
defined as follows by the action of the generators on the basis vectors $|\tau, k \rangle, k \in \mathbb{Z}$:
\begin{align} 
    J_0 |\tau, k \rangle &=(\tau +k) |\tau, k \rangle,  \label{IRzero} \\
    J_+ |\tau, k \rangle &= |\tau, k +1  \rangle,  \label{IRplus}\\
    J_- |\tau, k \rangle &= k(k - 1 + 2\tau) |\tau, k-1 \rangle. \label{IRminus}
\end{align}
(See in this connection \cite{howe2012non}.) It is readily checked that the Casimir element $J^2=J_0^2 - J_0 - J_+J_- = \tau(\tau - 1)$ on this representation space. The basis vectors are taken to be orthonormalized:
\begin{equation}
    \langle \tau,  k'| \tau,  k\rangle = \delta_{k'k}.
\end{equation}
\begin{rem}
Let us note the following.
\begin{enumerate}
    \item The representation defined above is not unitarisable \citep{tomasini2014unitary}.
    \item It is reducible and contains the unitary positive discrete series  \cite{vilenkin}, \cite{masson1991spectral}, \cite{groenevelt2001meixner}, \cite{howe2012non} as an irreducible component. This submodule is spanned by the basis vectors with $k \in \mathbb{Z}_+$.
\end{enumerate}
\end{rem}

Use now the formulas \eqref{embed_slL}, \eqref{embed_slM}, \eqref{embed_slX} of Proposition \ref{embedding} that define the isomorphism between $m\mathfrak{J}$ and $\mathfrak{su}(1,1)$ and take as already indicated $\tau=\frac{1}{2}(\alpha + \beta +1)$; the following actions of $L, M, X$ on the basis states $|\tau, k \rangle$ are readily found:
\begin{align}
  L |\tau, k \rangle=&\;(k + \alpha + 1) |\tau,  k \rangle \label{actL},  \\
  M |\tau,  k \rangle=& \;k \big[(k + \alpha + 1) |\tau,  k \rangle + ( k + \alpha + \beta) |\tau,  k - 1\rangle \big], \label{actM} \\
  X |\tau,  k \rangle=& \; |\tau,  k + 1 \rangle + |\tau,  k \rangle. \label{actX}
\end{align}

The adjoint actions can be read off directly:

\begin{align}
  L^T |\tau,  k \rangle=&\;(k + \alpha + 1) |\tau,  k \rangle \label{actLT},  \\ 
  M ^T|\tau,  k \rangle=& \;( k +1)(k + \alpha + \beta + 1) |\tau,  k +1 \rangle 
 + k(k + \alpha + 1 )|\tau,  k \rangle, \label{actMT} \\ 
  X^T |\tau, k \rangle=& \; |\tau,  k \rangle + |\tau,  k -1\rangle. \label{actXT}
\end{align}

Let us introduce the operator $\mathcal{T}_{\pm}$ on $\mathfrak{V}(\tau)$ such that:
\begin{equation}
    \mathcal{T}_{\pm}|\tau, k\rangle = |\tau, k \pm 1\rangle .
\end{equation}
Consider a vector $|f\rangle = \sum _{k=-\infty}^{\infty} f(k) |\tau, k \rangle$ in $\mathfrak{V}(\tau)$. We have
\begin{equation}
    \mathcal{T}_{\pm}|f\rangle = \sum _{k=-\infty}^{\infty} f(k) \mathcal{T}_{\pm}|\tau, k \pm 1\rangle = \sum _{k=-\infty}^{\infty} (T_{\mp}f(k))|\tau, k \pm 1\rangle,
\end{equation}
where $T_{\pm}$ stands for the shift operators acting on functions of $k$: $T_{\pm}f(k)=f(k\pm 1)$. 
\begin{rem} \label{3.2}
A realization of $m\mathfrak{J}$ in terms of shift operators can hence be inferred from the (dual) transformations of the components of a vector $|f\rangle$ in the basis $\{|\tau, k \rangle\}$ defined through: $V|f\rangle = \sum _{k=-\infty}^{\infty} f(k) V|\tau, k \rangle = \sum _{k=-\infty}^{\infty} (\mathrm{V}^Tf(k)) |\tau, k \rangle$. Equations \eqref{actLT}, \eqref{actMT}, \eqref{actXT} thus yield:
\begin{align}
\mathrm{L} &= (k +\alpha +1),  \label{shiftL} \\ 
\mathrm{M} &= ( k +1) ( k + \alpha + \beta + 1) T_+ + k(k + \alpha + 1), \label{shiftM} \\
\mathrm{X} &= T_- + 1 \label{shiftX}. 
\end{align}
The adjoints in this model are readily computed using $T_{\pm}^T=T_{\mp}$.
\end{rem}

We are now ready to construct the bases of $\mathfrak{V}(\tau)$ coming in adjoint pairs, whose overlaps will provide the algebraic interpretation we are looking for. (They will be in part the $d_n, d^*_n, e_z, e^*_z $ of the Introduction.)
The bases that will intervene are:
\begin{enumerate}
    \item The GEVP bases $\{|P_n\rangle\}$ and $\{|Q_n\rangle\}$:
    \begin{equation}
        M|P_n\rangle = \nu _n \;L|P_n\rangle \qquad M^T|Q_n\rangle = \nu _n \;L^T |Q_n\rangle. \label{GEVP}
    \end{equation}
    It will be recalled \cite{vinet2020unified}, \cite{zhedanov1999biorthogonal} that the sets $\{|P_n\rangle\}$ and $\{L^T|Q_n\rangle\}$ form by construction two biorthogonal ensembles of vectors:
    \begin{equation}
        \langle P_m| L^T|Q_n\rangle\ = 0,\quad  m\neq n.
    \end{equation}
    \item The EVP bases $\{|z\rangle\}$ and $\{\widetilde{|z\rangle}\}$:
        \begin{equation}
        X|z\rangle=z\,|z\rangle ,  \qquad X^T\widetilde{|z\rangle} =z\,\widetilde{|z\rangle}. \label{EVP}
    \end{equation}
    \item The EVP bases $\{|J_n\rangle\}$ and $\{\widetilde{|J_n\rangle}\}$:
    \begin{equation}
        M |J_n\rangle=\mu_n \,|J_n\rangle \qquad M^T\widetilde{|J_n\rangle}=\mu_n \, \widetilde{|J_n\rangle}. \label{EVPM}
    \end{equation}
\end{enumerate}



\subsection{Eigenvectors of $X$ and $X^T$}

It is directly checked that the EVP \eqref{EVP} are satisfied by
\begin{align}
    |z\rangle &= \gamma \sum_{k=-\infty}^{\infty} (z-1)^{-k-a}|\tau, k\rangle, \label{ez} \\
    \widetilde{|z\rangle} &= \Tilde{\gamma} \sum_{k=-\infty}^{\infty} (z-1)^{k+\Tilde{a}}|\tau, k\rangle, \label{ez*}
\end{align}
with $a, \Tilde{a} \in \mathbb{R}$ and where $\gamma, \Tilde{\gamma} \in \mathbb{C}$ are normalization constants. 
That $|z\rangle$ and $\widetilde{|z'\rangle}$ are orthogonal can be seen as follows. We have
\begin{equation}
    \widetilde{\langle z'|}z\rangle= \gamma \Tilde{\gamma}\sum_{k, l=-\infty}^{\infty} (z'-1)^{-k-a}(z-1)^{l+\Tilde{a}}\langle \tau, k|\tau, l\rangle. \label{last}
\end{equation}
Now let $z=1+e^{i\phi}$, $z'=1+e^{i\phi '}$, so that \eqref{last} becomes
\begin{equation}
    \widetilde{\langle z'|}z\rangle=\gamma \Tilde{\gamma}\;e^{i(\Tilde{a}\phi - a \phi ')}\sum_{k=-\infty} ^{\infty} e^{i(\phi - \phi ')k}.
\end{equation}
We then see that upon imposing
\begin{equation}
    a=\Tilde{a} + 1, \label{csts}
\end{equation}
we find 
\begin{equation}
     \widetilde{\langle z'|}z\rangle=-2\pi i \gamma \Tilde{\gamma} \delta(z-z') \label{delta}
\end{equation}
with the help of the Fourier series of Dirac's delta function and of a standard property of this distribution. Since $ \widetilde{\langle z'|}z\rangle$ is manifestly translation invariant, \eqref{delta} is preserved when the variable $z$ lies on the unit circle centered at $z=0$. 

We also have the completeness relation
\begin{equation}
    \frac{1}{2\pi i \gamma \Tilde{\gamma}} \oint _C dz \widetilde{|z \rangle} \langle z | = 1 \label{compl}
\end{equation}
where the contour $C$ consists in the unit circle infinitesimally deformed so that the singularity at $z=1$ lies inside $C$. Indeed, 
\begin{equation}
     \frac{1}{2\pi i \gamma \Tilde{\gamma}} \oint _C dz \widetilde{|z \rangle} \langle z |
     = \frac{1}{2\pi i} \oint _C dz (z-1)^{k-l-a +\Tilde{a}} \sum_{k, l = -\infty} ^{\infty}|\tau, k \rangle \langle \tau, l|.
\end{equation}
Again the choice \eqref{csts} for the integration constants $a$ and $\Tilde{a}$ consistently ensures that
the integral over $z$ becomes $\frac{1}{2\pi i} \oint _C dz (z-1)^{k-l-1} = \delta _{kl}$ and hence
\begin{equation}
     \frac{1}{2\pi i \gamma \Tilde{\gamma}} \oint _C dz \;\widetilde{|z \rangle}\langle z |= \sum _{k=-\infty}^{\infty} |\tau, k \rangle \langle \tau, k| =1.\label{complete}
\end{equation}
This will play a key role in the derivation of the orthogonality relations.

\subsection{GEVP bases}
We shall now obtain the bases $\{|P_n\rangle\}$ and $\{|Q_n\rangle\}$ of $\mathfrak{V}(\tau)$ that satisfy the GEVP \eqref{GEVP}. First we need to determine the set of eigenvalues $\nu$. From the explicit two-diagonal actions \eqref{actL}, \eqref{actM} of $L$ and $M$ on the basis vectors $\{|\tau, k\rangle\}$, it is readily seen that the (formal) determinantal
condition is
\begin{equation}
    det(M-\nu L) =\prod_{k=-\infty} ^{\infty} \left[ k(k + \alpha +1) - \nu \; ( k + \alpha + 1)\right] = 0
\end{equation}
and hence that the spectrum consists in the following values:
\begin{equation}
    \nu_n = n, \qquad n=0, \pm 1, \pm 2, \dots .
\end{equation}
\begin{rem}
In the following, as we consider GEVPs and EVPs, we shall limit ourselves to eigenvalues corresponding to \textit{non-negative} $n$, i.e. $n \in \mathbb{Z}^\geq$. This will not restrain the breadth of the algebraic description since the same results would be obtained with other choices. For completeness, indications on how the equations are handled for negative values of $n$ are given in Appendix \ref{neg}.
\end{rem}
Let 
\begin{equation}
    |P_n\rangle = \sum_{k=-\infty} ^{\infty} d_n(k) |\tau, k\rangle. \label{Pn}
\end{equation}
The generalized eigenvalue equation $M|P_n\rangle = n L |P_n\rangle $ implies the following recurrence relation for the expansion coefficients $d_n(k)$:
\begin{equation}
    (k + 1)(k + \alpha + \beta + 1)\; d_n(k+1) + (k - n) (k + \alpha + 1) \;d_n(k) = 0.\label{recurd}
\end{equation}
From \eqref{recurd}, it is immediately seen that for $n\geq0$,
\begin{equation}
    d_n(k) = 0 \quad \text{for} \quad k > n \quad \text{and} \quad k \in \mathbb{Z}_-.
\end{equation}
The explicit expression of the non-zero coefficients $d_n(k)$ reads
\begin{equation}
    d_n(k)= d_n(0) \frac{(-1)^k(-n)_k(\alpha +1)_k}{k!\;(\alpha + \beta +1)_k} \qquad  k=0, 1, 2, \dots, n.
      \label{sold}
\end{equation}

Turn now to the adjoint GEVP $M^T|Q_m\rangle = m \;L^T |Q_m\rangle$ which imposes on the coefficients $d_m^*(k)$ in
\begin{equation}
    |Q_m\rangle = \sum_{k=-\infty} ^{\infty} d_m^*(k) |\tau, k\rangle
\end{equation}
the recurrence relation
\begin{equation}
    k (k + \alpha + \beta) d_m^*(k - 1) + (k - m)( k + \alpha + 1) d_m^*(k) = 0. \label{recur*}
\end{equation}
Assuming as previously indicated, $m\geq 0$, one immediately notices that \eqref{recur*} implies 
\begin{equation}
    d_m^*(k) = 0 \quad \text{for} \quad k<m.
\end{equation}
In view of this fact, let
\begin{equation}
    k=l+m, \qquad l=0, 1, \dots,
\end{equation}
the relation \eqref{recur*} then becomes 
\begin{equation}
    (m + l) (l + m + \alpha + \beta) d_m^*(m+l-1) + l ( l + m + \alpha + 1) d_m^*(m+l). 
\end{equation}
It is found to have for solution
\begin{equation}
    d_m^*(m+l)= \frac{(-1)^l (m + 1)_l (m + \alpha + \beta + 1)_l}{l! (m + \alpha + 2)_l} \; d_m^*(m) \qquad l=0, 1, \dots . \label{d*}
\end{equation}
Apart from the initial condition $d_m^*(m)$, equation \eqref{d*} fully determines
\begin{equation}
|Q_m\rangle = \sum_{l=0} ^{\infty} d_m^*(m+l) |\tau, m + l\rangle.   \label{Qmf} 
\end{equation}

From general linear algebra considerations \cite{zhedanov1999biorthogonal}, \cite{tsujimoto2021algebraic}, \cite{vinet2020unified}, we know that the vectors $|P_n\rangle$ and $L^T|Q_m \rangle$ are biorthogonal for $n\neq m$. We have
\begin{align}
    &(\langle P_n|M)|Q_m\rangle =n (\langle P_n|L)|Q_m\rangle \nonumber \\
    &=\langle P_n|(M^T|Q_m\rangle)= m\langle P_m|(L^T|Q_m\rangle)= m (\langle P_n|L)|Q_m\rangle.
\end{align}
It follows that
\begin{equation}
    (n-m)(\langle P_n|L)|Q_m\rangle = (n-m)\langle P_n|(L^T|Q_m\rangle) = 0
\end{equation}
which implies the asserted biorthogonality if $m \neq n$.
Since the derivation we shall provide of the biorthogonality of the Askey polynomials will rest on this property, we shall next directly verify that it holds and determine the norm.

From the observations made above, we see that
\begin{align}
    \langle P_n|L^T|Q_m \rangle &= \sum_{k=-\infty}^n \sum_{l=0}^{\infty} d_n(k) d_m^*(l+m) \langle \tau, k|L^T|\tau, l+m\rangle \nonumber \\
    &=\sum_{k=-\infty}^n \sum_{l=0}^{\infty} d_n(k) d_m^*(l+m)(m + l + \alpha + 1)\; \delta _{k, l+m}. \label{PQ}
\end{align}
We readily find that
\begin{equation}
    \langle P_n|L^T|Q_m \rangle = 0 \qquad \text{if} \qquad m>n.
\end{equation}
It remains to consider the situation when $m \leq n$. Substituting in \eqref{PQ} the expressions \eqref{sold} and \eqref{d*} for $d_n(k)$  and $d_m^*(m+l)$, using a few properties of the Pochhammer symbols such as $x(x+1)_{l-1}=(x)_l$ and $(x)_{m+l}=(x)_m(x+m)_l$
and performing one of the sums, we arrive at
\begin{equation}
     \langle P_n|L^T|Q_m \rangle = d_n(0) d_m^*(m) \;(-1)^m \frac{(-n)_m (\alpha + 1)_{m+1} }{m!(\alpha + \beta + 1)_m } \sum_{l=0} ^{n-m} \frac{(-n+m)_l}{l!}.
\end{equation}
We then recall the following formula
\begin{equation}
    (1-x)^{\xi} = \sum _{k=0} ^{\infty} \frac{(-\xi)_k}{k!} x^k
\end{equation}
to conclude that
\begin{equation}
     \langle P_n|L^T|Q_m \rangle = 
     N_n\; \delta _{m,n}, \label{orth}
\end{equation}
with
\begin{equation}
    N_n = d_n(0) d_n^*(n) \;\frac{( \alpha + 1)_{n+1}}{(\alpha + \beta + 1)_n}. \label{N}
\end{equation}

\subsection{Askey polynomials and their biorthogonal partners}

Let us now identify some of the special functions that arise from this representation theoretic setting. In light of the completeness relation \eqref{compl} and the orthogonality relation \eqref{orth}, we see that $\widetilde{\langle z|} P_n \rangle$ and $\langle z|L^T|Q_n  \rangle$ provide two families of biorthogonal functions on the unit circle since 
\begin{equation}
    \frac{1}{2\pi i \gamma \Tilde{\gamma}} \oint _{|z|=1} dz \; \langle P_n \widetilde{| z\rangle}\langle z| L^T|Q_m \rangle = \langle P_n | L^T |Q_m\rangle = N_n\; \delta_{m,n}.\label{biorth}
\end{equation} 
These are explicitly obtained below.
\subsubsection{The overlaps $\widetilde{\langle z|} P_n \rangle$}
From  the expansions \eqref{ez*} and \eqref{Pn} of $ \widetilde{|z\rangle}$ and $|P_n \rangle$ over the orthonormal basis vectors $|\tau, k \rangle$ we have
\begin{equation}
    \widetilde{\langle z |}P_n \rangle = \Tilde{\gamma} \sum_{l=-\infty} ^n (z-1)^{l+\Tilde{a}}\; d_n(l).
\end{equation}
Upon inserting the expressions \eqref{sold} for $d_n(l)$, we observe that $\widetilde{\langle z |}P_n \rangle $ is the
${_2}F_1$ polynomial 
\begin{equation}
     \widetilde{\langle z |}P_n \rangle = \Tilde{\gamma}d_n(0) (z - 1)^{\Tilde{a}}\;  {_2}F_1\left({-n, \alpha + 1 \atop  \alpha + \beta + 1}; 1 - z \right ).
\end{equation}
The Askey polynomials are then recognized with the help of the following relation \cite{bateman1953higher}
\begin{equation}
 {_2}F_1\left({-n, \;b \atop  c}; z \right ) =  \frac{(c-b)_n}{(c)_n}{_2}F_1\left({-n, \;b \atop -n + b + 1 - c}; 1 - z \right ).   \label{Kummer1} 
\end{equation}
 We find
\begin{equation}
     \widetilde{\langle z |}P_n \rangle = \Tilde{\gamma} d_n(0) \frac{(\alpha +1)_n}{(\alpha + \beta +1)_n}(z - 1)^{\Tilde{a}} \;P_n(z; \alpha, \beta), \label{ztildePn}
\end{equation}
where the polynomials $P_n(z; \alpha, \beta)$ are as defined in \eqref{AP}. 
\begin{prop}
The Askey polynomials $P_n(z; \alpha, \beta)$ have a natural interpretation in the representation theory of the meta-Jacobi algebra. They occur according to \eqref{ztildePn} as the overlaps between two bases of the module $\mathfrak{V}(\tau =\frac{1}{2} (\alpha + \beta + 1))$ satisfying respectively equations defined in terms of the generators $X, L, M$ of $m\mathfrak{J}$. The first basis consists in the eigenvectors of $X^T$ (the transpose of $X$) and the second is formed by the vectors solving the GEVP defined by $L$ and $M$.
\end{prop}

\subsubsection{The overlaps $\langle z |L^T |Q_m \rangle$}
The biorthogonal partners to the Askey polynomials are obtained in a similar fashion. From \eqref{actLT}, \eqref{Qmf} and \eqref{d*} we have
\begin{equation}
    L^T |Q_m\rangle = d_m^*(m) (m+\alpha +1)\;\sum _{l=0}^{\infty} \frac{(-1)^l (m+1)_l (m+\alpha + \beta +1)_l}{l! (m+ \alpha + 1)_l} |\tau , l+m\rangle.
\end{equation}
Combining with \eqref{ez} and using the orthonormality of the basis vectors $|\tau, k\rangle$, we find
\begin{equation}
    \langle z |L^T | Q_m \rangle = \gamma d_m^*(m) (m + \alpha + 1) (z-1)^{-m-a}\; {_2}F_1\left({m+1,\; m+\alpha + \beta + 1 \atop m +\alpha + 1}; \frac{1}{1-z} \right ). \label{matQ}
\end{equation}
We may now use the fact that any three solutions of the hypergeometric equation are related by a linear relations and call upon transformation formulas of ${_2}F_1$ series under homographic transformations to make the biorthogonal partners of the Askey polynomials appear in this overlap. Indeed following the steps described in Appendix \ref{details1},
we arrive at the following expression: 
\begin{align}
    &\langle z| L^T| Q_m\rangle = \gamma d_m^*(m) (m + \alpha + 1) (z - 1)^{1-a}  \nonumber \\ 
    & \Bigg[
    \frac{\Gamma (m + \alpha +1 )\Gamma (\beta + 1)}{\Gamma (m + \beta + 2)\Gamma (\alpha)}\;{_2}F_1\left( {m + 1, \;1 - \alpha \atop m + \beta + 2}; z \right ) \nonumber \\
   & - \frac{\Gamma (m + \alpha + 1) \Gamma (m + \beta + 1)}{m! \; \Gamma (m + \alpha + \beta + 1)} \; (-z)^{-1-\beta} (1 - z)^{\alpha + \beta} Q_m (\frac{1}{z}, \alpha, \beta)  
   \Bigg] \label{overlap}
\end{align}
where $Q_m(z)$ is defined as in \eqref{part}.

\begin{rem}
Note that the first term in this expression for $\langle z| L^T| Q_m\rangle$ is a power series while the second one which contains the polynomial $Q_n$ in the variable $\frac{1}{z}$ has the transcendental factor  $z^{-\beta}$.
\end{rem}

Summing up:


\begin{prop} \label{partner}
The biorthogonal partners $Q_n(z, \alpha, \beta)$ of the Askey polynomials $P_n(z; \alpha, \beta)$ arise in the representation theory of the meta-Jacobi algebra in the overlaps, see \eqref{overlap}, between the eigenbasis vectors of the generatior $X$ and the basis vectors that obey the GEVP defined by the operators $M^T$ and $L^T$. 
\end{prop}

\subsubsection{Biorthogonality relation}
The interpretation of the Askey polynomials in the framework of the meta-Jacobi algebra leads to a natural derivation of their biorthogonality. 
Recall \eqref{biorth}. First observe that in multiplying the expressions of the overlaps $\widetilde{\langle z|}P_n \rangle$ and $\langle z | L^T |Q_m\rangle$,
as they are given by the formulas \eqref{ztildePn} and \eqref{overlap}, the factor $(z-1)^{1-a + \Tilde{a}}$ reduces to $1$ because of \eqref{csts}. Furthermore, one hence observes that the product of the first term in \eqref{overlap} - a power series - with the polynomial $P_n(z, \alpha, \beta)$ will give a vanishing contribution when integrated over the circle $|z|=1$. Equation \eqref{biorth} thus yields
\begin{align}
    d_n(0) d_m^*(m) &\frac{(m + \alpha + 1) (\alpha + 1)_n}{(\alpha + \beta + 1)_n} \frac{\Gamma (m + \alpha + 1) \Gamma (m + \beta + 1)}{m!\; \Gamma (m + \alpha + \beta + 1)} \nonumber\\
    &\times \frac{-1}{2 \pi i} \oint_{|z|=1} dz \;
    (-z)^{-1 - \beta} (1-z)^{\alpha + \beta}P_n(z, \alpha, \beta) Q_m (\frac{1}{z}, \alpha, \beta) = N_n \delta_{m,n}.
\end{align}
   Mindful of formula \eqref{N} for $N_n$, we thus recover precisely the biorthogonality relation \eqref{bi}.

\subsection{Eigenbases of $M$ and $M^T$}
We now undertake to show that the Jacobi polynomials can be described within the same algebraic framework. We already noted in Proposition \ref{embedding} that the elements $X$ and $M$ generate the Jacobi algebra $\mathfrak{J}$ which is thus embedded in $m\mathfrak{J}$. We therefore expect to see the Jacobi polynomials occur in the overlaps between the eigenvectors of $X$, $X^T$ and $M^T$, $M$ respectively. We shall hence first determine the eigenbases of  $\mathfrak{V}(\tau)$ associated to $M$ and $M^T$.

From the two-diagonal action \eqref{actM} of $M$, we see that the spectrum $\{\mu _n \}$ of this operator is of the form
\begin{equation}
    \mu_n = n(n + \alpha + 1).
\end{equation}
Consider the EVPs \eqref{EVPM}. Set,
\begin{equation}
    |J_n\rangle = \sum_{k=-\infty}^{\infty} f_n (k) |\tau, k\rangle. \label{Jn}
\end{equation}
The eigenvalue equation $M |J_n\rangle = n(n + \alpha + 1) |J_n\rangle$ yields the following two-term recurrence relation for the coefficients $f_n(k)$:
\begin{equation}
    (k + 1 + \alpha + \beta) f_n (k + 1) + \left[ k(k + \alpha + 1) - n(n + \alpha + 1) \right] f_n(k) = 0
\end{equation}
which can be rewritten as
\begin{equation}
    (k + 1) (k + 1 + \alpha + \beta) f_n (k + 1) + (k - n) (k + n + \alpha + 1) f_n (k) = 0.
\end{equation}
Here again, we shall focus on the case $n \geq 0$. The above equation is then found to imply that 
\begin{equation}
    f_n(k) = 0  \quad \text{for} \quad k>n \quad \text{and} \quad k \in \mathbb{Z}_-
\end{equation}
and is solved by
\begin{equation}
    f_n(k) = f_n (0) (-1)^k\frac{(-n)_k (n + \alpha + 1)_k}{k! (\alpha + \beta + 1)_k}, \qquad k = 0, 1, \dots, n. \label{fn}
\end{equation}
Similarly, let
\begin{equation}
    \widetilde{|J_n\rangle} = \sum_{k=-\infty}^{\infty} \Tilde{f}_n (k) |\tau, k\rangle. \label{Jnt}
\end{equation}
The EVP $M^T \widetilde{|J_n \rangle} = n (n + \alpha + 1) \widetilde{|J_n \rangle} $ is readily seen to give:
\begin{equation}
    k ( k + \alpha + \beta) \Tilde{f}_n (k - 1) + (k - n) (k + n + \alpha + 1) \Tilde{f}_n (k)= 0. \label{recurftilde}
\end{equation}
In this instance, for $m \geq 0$, we observe that
\begin{equation}
    \Tilde{f}_n (k) = 0 \qquad \text{when } \qquad k<n.
\end{equation}
We set $k = n + l , l = 0, 1, 2, \dots$ and convert \eqref{recurftilde} into
\begin{equation}
    (n + l) (n + l + \alpha + \beta) \Tilde{f}_n (n + l - 1) + l ( l + 2n + \alpha + 1) \Tilde{f}_n (n + l) = 0
\end{equation}
to find 
\begin{equation}
    \Tilde{f}_n (n + l) = \Tilde{f}_n (n) (-1)^l \frac{(n + 1)_l (n +1 + \alpha + \beta)_l}{l! (2n + \alpha + 2)_l}. \label{fnt}
\end{equation} 

We may now verify that $|J_n\rangle$ and $\widetilde{|J_m \rangle}$ are orthogonal when $m \neq n$. This proceeds in a way similar to the computation of $\langle P_n | L^T | Q_m \rangle $ carried out before. Clearly, $\langle J_n \widetilde{| J_m \rangle} = 0$ if $m > n$. If $m \leq n$, after some algebraic simplifications, we see that
\begin{equation}
    \langle J_n \widetilde{| J_m \rangle} = f_n(0) \Tilde{f}_m(n) (-1)^m\frac{(-n)_m (n + \alpha + 1)_m}{m! (\alpha + \beta + 1)_m }   {_2}F_1\left({m - n, \;n + m + \alpha + 1 \atop 2m + \alpha + 2}; 1 \right ).
\end{equation}
From the Vandermonde formula \cite{gasper2004basic} 
\begin{equation}
      {_2}F_1\left({-n, \;b \atop c}; 1 \right ) = \frac{(c-b)_n}{(c)_n}
\end{equation}
we may then conclude that owing to the factor $(m - n +1)_{n-m}$ that appears,  $\langle J_n \widetilde{| J_m \rangle} = 0$ unless $n=m$, in which case
\begin{equation}
     \langle J_n \widetilde{| J_m \rangle} = \mathcal{N}_n \delta_{m,n},
\end{equation}
with
\begin{equation}
    \mathcal{N}_n = f_n(0) \Tilde{f}_n(n) \frac{(n + \alpha + 1)_n}{(\alpha + \beta + 1)_n}.\label{normJ}
\end{equation}
\subsection{Jacobi polynomials}
We will now observe how the Jacobi polynomials emerge in this framework and indicate how this allows for another derivation of their orthogonality relation.
\subsubsection{The overlaps}
Let us now look at the overlaps $\widetilde{\langle z |} J_n \rangle$ and $\langle z \widetilde{| J_m \rangle}$. From \eqref{ez*}, \eqref{Jn} and \eqref{fn} we obtain
\begin{equation}
   \widetilde{\langle z |} J_n \rangle = \Tilde{\gamma} f_n(0) (z - 1 )^{\Tilde{a}} \;
   {_2}F_1\left({-n, \;n + \alpha + 1 \atop \alpha + \beta + 1}; 1-z \right ). \label{over}
\end{equation}
Using \eqref{comeon} (or equivalently \eqref{Kummer1})
we arrive at 
\begin{equation}
    \widetilde{\langle z |} J_n \rangle = \Tilde{\gamma} (z - 1)^{\Tilde{a}} \frac{(-1)^n \;\Gamma (1 + \alpha + \beta) \Gamma (\beta) \Gamma (2n + \alpha + 1) \Gamma (1 - \beta)}{\Gamma (-n + \beta) \Gamma (n + \alpha + \beta + 1)\Gamma (n + \alpha + 1) \Gamma (n + 1 - \beta)}  \;\hat{P}_n^{(\alpha, \beta)}(z), \label{overJ}
\end{equation}
where $\hat{P}_n^{(\alpha, \beta)}(z)$ are the Jacobi polynomials defined in \eqref{Jac} extended to the complex plane.

The second overlap is recovered from \eqref{ez}, \eqref{Jnt} and  \eqref{fnt}. We find
\begin{equation}
    \langle z \widetilde{| J_m \rangle} = \gamma \Tilde{f}_m(m) (z-1)^{-m - a} \; {_2}F_1\left({m + 1, \;m + \alpha + \beta + 1 \atop 2m +\alpha + 2}; \frac{1}{1-z} \right ). \label{secondoverlap}
\end{equation}
At this point by performing the transformations described in Appendix \ref{compJac} that make use of identities involving gamma functions and solutions of the hypergeometric equation, the following formula is discovered:
\begin{align}
    & \langle z \widetilde{| J_m \rangle} =  (-1)^{m + 1}\gamma \Tilde{f}_m(m) (z-1)^{1 - a}\nonumber \\ 
    & \Bigg[
    \frac{\Gamma (2m + \alpha +2 )\Gamma (-\beta)}{\Gamma (m + \alpha + 1)\Gamma (m -\beta +1)} \;{_2}F_1\left( {m + 1, \;-m - \alpha \atop 1 + \beta}; z \right ) + \nonumber \\
   & \frac{(-1)^m \Gamma (2m + \alpha + 2) \Gamma (\beta) \Gamma (2m + \alpha + 1)\Gamma (1 - \beta)}{m! \; \Gamma (m + \alpha + \beta + 1) \Gamma (m + \alpha + 1) \Gamma (m + 1 - \beta)} \; (-z)^{-\beta} (1 - z)^{\alpha + \beta} \; \hat{P}_m^{(\alpha, \beta)}(z) 
   \Bigg]. \label{overlap2}
\end{align}

\begin{rem}
As already encountered in the expression \eqref{overlap} of $\langle z| L^T| Q_m\rangle$, we see that the first term in \eqref{overlap2} is a power series while the second that involves the Jacobi polynomial contains the transcendental term $z^{-\beta}$.
\end{rem}
\begin{prop}
The Jacobi polynomials $\hat{P}_m^{(\alpha, \beta)}(z)$ over $\mathbb{C}$ also arise in the context of the meta-Jacobi algebra $m\frak{J}$. They occur as per the equations  \eqref{overJ} and \eqref{overlap2} in two overlaps between eigenbases of the module $\mathfrak{V}(\frac{1}{2} (\alpha + \beta + 1))$ : on the one hand between the eigenstates of $M$ and $X^T$ and on the other hand between those of $M^T$ and $X$.
\end{prop}

\subsubsection{Orthogonality}
This interpretation of the Jacobi polynomials in the framework of the algebra $m\mathfrak{J}$ entails a derivation of their orthogonality. Owing to the completeness relation \eqref{complete}, we have
\begin{equation}
    \frac{1}{2 \pi i \gamma \Tilde{\gamma}} \oint _{C_{|z|=1} }\langle J_n \widetilde{| z \rangle} \langle z \widetilde{| J_m \rangle} = \langle J_n \widetilde{| J_m \rangle} = \mathcal{N}_n \delta_{m,n} \label{perpJ}
\end{equation}
where $\mathcal{N}_n$ is given by \eqref{normJ}. 
When substituting the expressions \eqref{overJ} and \eqref{overlap2} for $\widetilde{\langle z |} J_n \rangle$ and $\langle z \widetilde{| J_m \rangle}$, we first observe anew that the resulting factor $(z-1)^{1-a + \Tilde{a}} = 1$ since $1 - a + \Tilde{a} = 0$. Then, we note that the product of $\widetilde{\langle z |} J_n \rangle$ with the first term of \eqref{overlap2} is a power series that will integrate to zero over the unit circle. Taking into account the formula \eqref{normJ} for $\mathcal{N}_n$ and after some simplifications, equation \eqref{perpJ} thus amounts to 
\begin{equation}
    - \frac{1}{2 \pi i} \Big(\frac{\pi}{\sin \pi \beta}\Big) \oint _{C_{|z|=1}} \;dz (-z)^{-\beta} (1 - z) ^{\alpha + \beta} \hat{P}_n^{(\alpha, \beta)}(z)\hat{P}_m^{(\alpha, \beta)}(z) = h_n \; \delta _{m,n}, \label{quasiorth}
\end{equation}
with $h_n$ given in \eqref{hn}. In obtaining \eqref{quasiorth} we have used the identity $\Gamma(x) \Gamma (1-x) = \frac{\pi}{\sin \pi x}$ and in particular
\begin{equation}
    \Gamma (-n + \beta) \Gamma (n + 1 - \beta) = \frac{\pi}{\sin \pi (-n + \beta)} = (-1)^n \dfrac{\pi}{\sin \pi \beta}.
\end{equation}

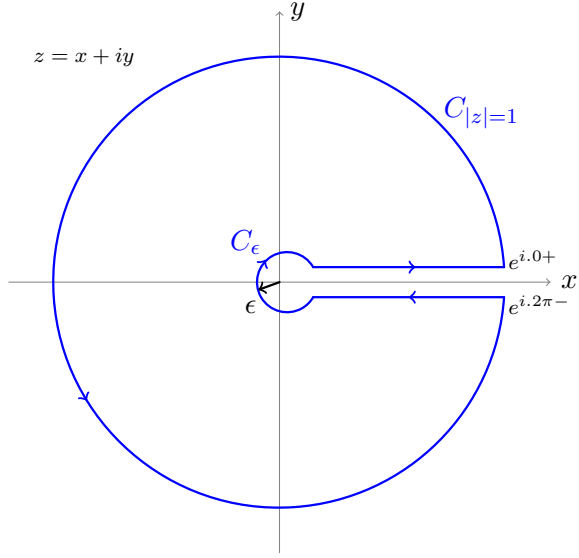
\begin{figure}
   \centering
 
\begin{tikzpicture}[scale=2,
                   line width=0.3mm,
                    line cap=round,
                    dec/.style args={#1#2}{
                   decoration={markings, mark=at position #1 with {#2}},
                    postaction={decorate}                                   }]
 Axes
\path[gray,very thin] (-1.8,0) edge[->] (1.8,0);
\path[gray,very thin] (0,-1.8) edge[->] (0,1.8);
\draw (1.8, 0) node[right]{$x$};
\draw (0, 1.8) node[right]{$y$};
\draw (-1.3,1.5) node{\scriptsize $z=x+iy$};
 Contours
Big circle, radius=1
\draw[blue,dec={0.59}{\arrow{>}}]  
    (1.5*0.99500416527,0.09983341664)coordinate(1a) arc ({asin(0.06666666666)}:{360-asin(0.06666666666)}:1.5cm)coordinate(2a); 
 Small circle, radius=0.2
\draw[blue,dec={0.39}{\arrow{<}}]  
    (0.22353234609,0.09983341664)coordinate(3a) arc
    ({asin(0.5)}:{360-asin(0.5)}:0.2cm)coordinate(4a); 
\draw (40: 1.75) node{\BLUE{$C_{|z|=1}$}};
\draw (130: 0.35) node{\BLUE{$C_{\epsilon}$}};
\draw[blue, dec={0.5}{\arrow{<}}] (1a)--(3a);
\draw[blue, dec={0.5}{\arrow{>}}] (2a)--(4a);
\draw (1.45,0.15) node[right]{\scriptsize $e^{i.0+}$};
\draw (1.45,-0.15) node[right]{\scriptsize $e^{i.2\pi-}$};
 Radius
\path (0,0) edge[->] (200:0.15);
\draw (220:0.25) node{$\epsilon$};
\end{tikzpicture}
\caption{The contour $\Xi$}
   \label{fig:contour}
\end{figure}

Finally, the orthogonality of the Jacobi polynomials on the interval $[0, 1]$ is recovered by employing the contour depicted in Figure \ref{fig:contour} and computations carried out in \cite{hendriksen1986orthogonal}. Schematically the contour $\Xi = C_{|z|=1} + [1, 0] + C_{\epsilon} + [0, 1]$, it is composed of the unit circle (short of crossing the branch cut), the segment from $x=1$ to $x=0$ below the branch cut, a circle of radius $\epsilon$ around $z=0$ and the segment from $x=0$ to $x=1$ above the branch cut. Consider the integral in \eqref{quasiorth} with the contour $C_{|z|=1}$ replaced by the contour $\Xi$ of figure \ref{fig:contour}. Since no singularities are enclosed by $\Xi$ that integral is equal to zero.

If we restrict $\beta$ to be smaller than $1$, i.e. if we take $\beta < 1$ as in the standard definition of the Jacobi polynomials, it is readily seen that
\begin{equation}
    \lim_{\epsilon \to 0} \frac{1}{2 \pi i} \oint _{C_{\epsilon}}\;dz (-z)^{-\beta} (1 - z) ^{\alpha + \beta} \hat{P}_n^{(\alpha, \beta)}(z)\hat{P}_m^{(\alpha, \beta)}(z) = 0, \qquad \text{for} \;\beta < 1.
\end{equation}
It follows that the integral over $C_{|z|=1}$ must be the negative of the sum of the integrals over and above the real axis.
Hence, recalling the choice of branch $(-z)^{-\beta} = |z|^{-\beta}$ when $\arg z = \pi$, we have
\begin{align}
    &- \frac{1}{2 \pi i} \Big(\frac{\pi}{\sin \pi \beta}\Big) \oint _{C_{|z|=1}} \;dz (-z)^{-\beta} (1 - z) ^{\alpha + \beta} \hat{P}_n^{(\alpha, \beta)}(z)\hat{P}_m^{(\alpha, \beta)}(z) = \nonumber \\
    & \frac{1}{2i} \Big(\frac{1}{\sin \pi \beta}\Big) \; e^{i\pi \beta} (1 - e^{-2 \pi i \beta)} \int_0^1\; x^{-\beta} (1 - x) ^{\alpha + \beta} \hat{P}_n^{(\alpha, \beta)}(x)\hat{P}_m^{(\alpha, \beta)}(x)\; dx.
\end{align}
The factors before the integral sign in the last expression cancel and this gives the orthogonality relation \eqref{orthrelJac} of the Jacobi polynomials in view of \eqref{quasiorth}
 
\section{Algebraic derivation of the properties of the Askey polynomials}\label{4}

We shall indicate in this section how various properties of the biorthogonal Askey polynomials on the circle naturally follow from their interpretation based on the meta-Jacobi algebra.
Recall that
\begin{equation}
    |P_n\rangle = d_n(0) \sum_{k=0}^n (-1)^k \frac{(-n)_k (\alpha + 1)_k}{k! (\alpha + \beta + 1)_k} \;|\tau, k\rangle. \label{P}
\end{equation}
Looking at the overlap $\widetilde{\langle z|} P_n\rangle$ given in \eqref{ztildePn}, without loss of generality, we can set from now on:
\begin{equation}
    \Tilde{\gamma} = 1, \qquad \Tilde{a} = 0, \qquad a=1. \label{choices}
\end{equation}
It is moreover natural to take the initial values $d_n(0)$ of the recurrence relation \eqref{recurd} to be
\begin{equation}
    d_n(0; \alpha, \beta) = \frac{(\alpha + \beta + 1)_n}{(\alpha + 1)_n}
\end{equation}
so that
\begin{equation}
    \widetilde{\langle z|} P_n\rangle = P_n(z; \alpha, \beta), \label{Askpol}
\end{equation}
identifying $\widetilde{\langle z|} P_n\rangle$ precisely with the Askey polynomials. This also means that $|\tau, n\rangle$ has coefficient $1$ in $|P_n\rangle$:
\begin{equation}
    |P_n\rangle = \frac{(\alpha + \beta + 1)_n}{(\alpha + 1)_n} |\tau, 0\rangle + \dots + |\tau, n\rangle.
\end{equation}
\subsection{Action of $L$ and $R=XL$ in the basis $\{|P_n\rangle \}$}
We shall now show that the generator $L$ and the product $XL$ act in a two-diagonal fashion in the basis $\{|P_n\rangle, n=0, 1, \dots\}$. We have
\begin{equation}
    L |P_n(\alpha, \beta)\rangle = \frac{(\alpha + \beta + 1)_n}{(\alpha + 1)_n} \sum_{k=0}^n (-1)^k \frac{(-n)_k (\alpha + 1)_k}{k! (\alpha + \beta + 1)_k} L|\tau, k\rangle.
\end{equation}
From eq.\eqref{actL} which reads $L |\tau, k\rangle = (k + \alpha + 1) |\tau, k\rangle$, and the identity
\begin{align}
    (-n)_k (k + \alpha + 1) &= (-n)_k \big[n + \alpha + 1 + (-n + k) \big] \\ \nonumber
    &= (n + \alpha + 1) (-n)_k - n(-n + 1)_k,
\end{align}
we see that 
\begin{equation}
    L |P_n(\alpha, \beta)\rangle = (n + \alpha + 1) |P_n(\alpha, \beta)\rangle - \frac{n (n + \alpha + \beta)}{(n + \alpha)} |P_{n-1} (\alpha, \beta) \rangle. \label{LonP}
\end{equation}

Alternatively, using
\begin{equation}
    (k + \alpha + 1) (\alpha + 1)_k = (\alpha + 1)_{k+1} = \big((\alpha +1) + 1\big)_k (\alpha + 1), 
\end{equation}
we note that $L$ has also the effect of shifting the parameters:
\begin{equation}
    L |P_n(\alpha, \beta)\rangle = ( n + \alpha + 1)\; |P_n(\alpha + 1, \beta -1)\rangle. \label{contigL}
\end{equation}

Consider now the action of the operator $R=XL$. Knowing that $L$ acts diagonally as per \eqref{actL} on the basis vectors $|\tau, k\rangle$ and that according to \eqref{actX}, i.e. $X|\tau , k\rangle = |\tau , k\rangle + |\tau , k + 1\rangle$, we have
\begin{equation}
        R |P_n(\alpha, \beta)\rangle = \frac{(\alpha + \beta + 1)_n}{(\alpha + 1)_n} \sum_{k=0}^n (-1)^k \frac{(-n)_k (\alpha + 1)_k (k + \alpha + 1)}{k! (\alpha + \beta + 1)_k} \Big[|\tau, k\rangle + |\tau, k + 1\rangle\Big].
\end{equation}
Collecting the factors of the vectors $|\tau, k\rangle, \; k=0, \dots, n+1$, we find
\begin{align}
   & R |P_n(\alpha, \beta)\rangle = \nonumber \\ \nonumber
   &\frac{(\alpha + \beta + 1)_n}{(\alpha + 1)_n} \Bigg[ (\alpha + 1) |\tau, 0\rangle \\ 
   &+ \sum_{k=1}^n  \frac{(-1)^k (\alpha + 1)_k}{k! (\alpha + \beta + 1)_k} \Big[ (-n)_k (k + \alpha + 1) - (-n)_{k-1} k(k + \alpha + \beta) \Big] \; |\tau,k \rangle \label{XLP}\\ \nonumber
   & + \frac{(\alpha + 1)_{n+1}}{(\alpha + \beta + 1)_n} \; |\tau, n + 1 \rangle \Bigg]. 
\end{align}
The two relations
\begin{align}
    &k (-n)_{k-1} = (-n)_k - (-n-1)_k, \label{1*}\\
    &(-n)_k (-n-1) - (-n - 1)_k (n + 1 -k)=0 \label{2*}
\end{align}
come in handy in deriving the following identity:
\begin{align}
    &(-n)_k (k + \alpha + 1) - k (-n)_{k-1} (k + \alpha + \beta)  \nonumber \\ \nonumber
    & = (-n)_k (k + \alpha + 1) - \big[ (-n)_k - (-n-1)_k \big] (k + \alpha + \beta)  \nonumber \\
    & = (-n)_k (1 - \beta) + (-n - 1)_k (k + \alpha + \beta) \nonumber \\
    & = (-n)_k (-n - \beta) + (-n-1)_k ( n + \alpha + \beta + 1).\label{line4}
\end{align}
Clearly, \eqref{1*} has been used in getting the second line and \eqref{2*} has been added to the third line to obtain the end result. Upon inserting this relation \eqref{line4} in \eqref{XLP}, we recognize easily that $R$ is a two-diagonal raising operator:
\begin{equation}
    R |P_n(\alpha, \beta)\rangle = (n + \alpha + 1) \;|P_{n+1}(\alpha, \beta)\rangle - (\beta + n) \;|P_n(\alpha, \beta)\rangle. \label{2dR}
\end{equation}

In the following we shall also consider the element
\begin{equation}
    \Tilde{R} = XM.
\end{equation}
\begin{rem} \label{4.1}
Given that the vectors $|P_n(\alpha, \beta)\rangle $ satisfy the GEVP $\;M|P_n(\alpha, \beta)\rangle = n L |P_n(\alpha, \beta)\rangle$, the eqs. \eqref{LonP} and \eqref{2dR} readily provide the actions of $M$ and $\Tilde{R}$ on these vectors.
\end{rem}

\subsection{A differential realization}

A differential model of the meta-Jacobi algebra is directly obtained. With the choices \eqref{choices} we have
\begin{equation}
    \widetilde{\langle z |} \tau, k\rangle \equiv f(z, k) = (z - 1)^k.
\end{equation}
We can dually define an operator acting on the variable $z$ as follows
\begin{equation}
    \mathcal{O}_z \widetilde{\langle z |} \tau, k\rangle = \widetilde{\langle z |} O |\tau, k\rangle = \mathrm{O}_k \;f(z, k),
\end{equation}
where $\mathcal{O}_z$ corresponds to the operator $O$ acting on the module $\mathfrak{V}(\frac{1}{2} (\alpha + \beta + 1))$ and $\mathrm{O}_k$ as in remark \ref{3.2},
acts on the components $f(z,k)$ of the vector $\widetilde{|z\rangle}$. With $O = L, M, X$ we find:
\begin{prop}
The differential operators $\mathcal{L}$, $\mathcal{M}$ and $\mathcal{X}$ given below provide a realization of the commutation relations \eqref{com_LM}, \eqref{com_LX} and \eqref{com_MX} of the meta-Jacobi algebra.
\begin{equation}
 \mathcal{L} = (z-1) \partial_z + (\alpha+1) \mathcal{I}; \label{L1}   
\end{equation}
\begin{equation}
\mathcal{M}= z(z-1) \partial_z^2 + \left[ (\alpha+2)z + \beta-1   \right] \partial_z; \label{L2} 
\end{equation}
\begin{equation}
    \mathcal{X}=z. \label{z}
\end{equation}

\end{prop}
It follows that $R=XL$ and $\Tilde{R}=XM$ are realized by
\begin{equation}
    \mathcal{R} = z(z-1)\partial _z + (\alpha + 1) z; \label{R1}
\end{equation}
\begin{equation}
    \Tilde{\mathcal{R}} = z^2 (z-1) \partial _z ^2 + z\left[(\alpha+2)z+\beta-1 \right] \partial_z \label{R2} .
\end{equation}
\begin{rem} \label{after}

One may also take $X$ acting on the left on $\widetilde{\langle z |}$ and giving the eigenvalue $z$ so that
\begin{equation}
    \mathcal{R} \;\widetilde{\langle z |}P_n(\alpha, \beta)\rangle = \widetilde{\langle z |}XL |P_n(\alpha, \beta)\rangle= z \widetilde{\langle z |}L |P_n(\alpha, \beta)\rangle  = z\mathcal{L} \;\widetilde{\langle z |}P_n(\alpha, \beta)\rangle \label{RzL} 
\end{equation}
and similarly for $\Tilde{R} = XM$.

\end{rem}
\begin{rem}
Note that $\mathcal{L}$ and $\mathcal{M}$ have the property of stabilizing spaces of polynomials of given degrees while $\mathcal{X}$, $\mathcal{R}$, $\Tilde{\mathcal{R}}$ raise the degree by $1$. In the spirit of studies carried out in \cite{grunbaum2017tridiagonalization}, \cite{tsujimoto2019rational}, \cite{gaboriaud2020degenerate} for example, $\mathcal{X}$, $\mathcal{R}$, $\Tilde{\mathcal{R}}$ are operators of Heun type.
\end{rem}
\begin{rem}
Observe that $\mathcal{M}$ precisely coincides with the hypergeometric operator \eqref{M} albeit in the variable $z$.
\end{rem}
\begin{rem}
This differential model for $m\mathfrak{J}$ can also be retrieved by using the Barut-Ghirardello (B-G) realization of $\mathfrak{su}(1,1)$:
\begin{align}
    J_0=&\;(z-1)\frac{d}{dz} + \tau, \nonumber \\
    J_+=&\;(z-1),\label{sl2real} \\
    J_-=&\;(z-1)\frac{d^2}{dz^2}+2\tau\frac{d}{dz}, \qquad \tau = \frac{1}{2}(\alpha + \beta + 1),\nonumber
\end{align}
in the formulas \eqref{embed_slL}, \eqref{embed_slM}, \eqref{embed_slX} giving $L$, $M$, and $X$ in terms of the $\mathfrak{su}(1,1)$ generators. Note that the variable $z$ is here translated by $1$ with respect to the usual B-G formulas.
\end{rem}
Given the expression \eqref{Askpol} of $P_n(z; \alpha, \beta)$ as $\widetilde{\langle z |}P_n(\alpha, \beta) \rangle$, in view of the actions \eqref{LonP}, \eqref{2dR} of $L$ and $R$ on the vectors $|P_n\rangle$, of remark \ref{4.1} and of the realizations given above ( \eqref{L1}, \eqref{L2}, \eqref{z} \eqref{R1}, \eqref{R2}) of these operators, we have the following.
\begin{prop}
The biorthogonal Askey polynomials $P_n(z; \alpha, \beta)$ on the unit circle satisfy the following differential identities:
\begin{align}
\mathcal{L} P_n(z; \alpha, \beta) &= (n+\alpha+1) P_n(z; \alpha, \beta) - \frac{n (\alpha+\beta+n)}{\alpha+n} P_{n-1}(z; \alpha, \beta), \label{L1_rec} \\
\mathcal{M} P_n(z; \alpha, \beta) &= n(n+\alpha+1) P_n(z; \alpha, \beta) - \frac{n^2 (\alpha+\beta+n)}{\alpha+n} P_{n-1}(z; \alpha, \beta), \label{L2_rec} \\
\mathcal{R} P_n(z; \alpha, \beta) &= (n+\alpha+1) P_{n+1}(z; \alpha, \beta) - (\beta+n) P_n(z; \alpha, \beta), \label{R1_rec} \\
\Tilde{\mathcal{R}} P_n(z; \alpha, \beta) &= n(n+\alpha+1) P_{n + 1}(z; \alpha, \beta) - n(\beta+n) P_n(z; \alpha, \beta) \label{R2_rec} .
\end{align}
\end{prop}

\subsection{Bispectrality}
The bispectral equations of the Askey polynomials can now easily be identified and interpreted in terms of generalized eigenvalue problems.
\subsubsection{The differential equation}
The GEVP $M |P_n(\alpha, \beta)\rangle = n L |P_n(\alpha, \beta)\rangle$ translates after projection on $\widetilde{\langle z |}$ into the second order differential equation:
\begin{equation}
    \mathcal{M} P_n(z; \alpha, \beta) = n \mathcal{L} P_n(z; \alpha, \beta) \label{diffeq}
\end{equation}
with eigenvalue $n$ and where the operators $\mathcal{M}$ and $\mathcal{L}$ are respectively given by \eqref{L2_rec} and \eqref{L1_rec}.

\subsubsection{The recurrence relation}
The recurrence relation is obtained by considering the GEVP
\begin{equation}
    \mathcal{R} P_n(z; \alpha, \beta) = z \mathcal{L} P_n(z; \alpha, \beta) \label{GEVPRL}
\end{equation}
which is satisfied by construction (see \eqref{RzL} in Remark \ref{after}). Expressing in \eqref{GEVPRL} the two- diagonal actions \eqref{L1_rec} and \eqref{R1_rec} of $\mathcal{L}$ and $\mathcal{R}$, one arrives at the recurrence relation
\begin{equation}
P_{n+1}(x) + b_n P_n(x) = x \left(P_n(x) + g_n P_{n-1}(x)  \right) \label{rec_P} 
\end{equation} 
where
\begin{equation}
b_n = -\frac{\beta+n}{\alpha+n+1}, \quad g_n = -\frac{n(n+\alpha+\beta)}{(\alpha+n)(\alpha+n+1)}. \label{bg_rec} 
\end{equation}
This recurrence relation was obtained by Hendriksen and van Rossum in \cite{hendriksen1986orthogonal}. It was derived in \cite{grunbaum2004linear} by considering linear pencils in 
$\mathfrak{sl}_2$. It is also constructed through a gluing procedure by Kim and Stanton in their recent study of $R_I$ polynomials \cite{kim2020orthogonal}. 
\begin{rem}
It is manifest from this recurrence relation of $R_I$ - type \cite{ismail1995generalized}, that $z$ (resp. $X$) is a lower Hessenberg matrix on the space of Askey polynomials (resp. in the basis $\{|P_n(\alpha, \beta)\rangle \}$). This feature of the representation theory of $m\mathfrak{J}$ was also observed in the study of the meta-Hahn algebra \cite{vinet2020unified}.
\end{rem}

\begin{prop}
The biorthogonal Askey polynomials defined on the unit circle are bispectral. They satisfy the differential equation \eqref{diffeq} and the recurrence relation of $R_I$ - type \eqref{rec_P} with coefficients \eqref{bg_rec}. Both spectral equations are of GEVP type.
\end{prop}

\subsection{Contiguity relations}
Some contiguity relations for the Askey polynomials arise also naturally in the meta-Jacobi algebra framework. Indeed, we already observed in eq. \eqref{contigL} that the generator $L$ has the effect of performing the shifts $\alpha \rightarrow \alpha + 1$, $\beta \rightarrow \beta - 1$ when acting on $|P_n (\alpha, \beta) \rangle$. That $M$ has a similar effect follows from the fact that $M=nL$ in the GEVP basis $|P_n (\alpha, \beta) \rangle$. This translates into the following for the polynomials $P_n(\alpha, \beta) = \widetilde{\langle z |} P_n (\alpha, \beta)\rangle$.
\begin{prop}
The Askey polynomials $P_n(\alpha, \beta)$ verify the following contiguity equations:
\begin{align}
\mathcal{L} P_n(z;\alpha,\beta) &= (\alpha+n+1) P_n(z; \alpha+1, \beta-1),  \label{Lalphab} \\
\mathcal{M} P_n(z;\alpha,\beta) &= n(\alpha+n+1) P_n(z; \alpha+1, \beta-1),\label{Malpha}
\end{align}
where $\mathcal{L}$ and $\mathcal{M}$ are the differential operators \eqref{L1} and \eqref{L2} respectively.
\end{prop}
\begin{rem}
Given the explicit form \eqref{AP} of the Askey polynomials, the above relations can be checked directly on $P_n(z;\alpha,\beta)$ with the differential operators $\mathcal{L}$ and $\mathcal{M}$. Having done this, comparing \eqref{Lalphab} and \eqref{Malpha} offers a way to show that the Askey polynomials are solutions of the GEVP \eqref{diffeq}.
\end{rem}

\subsection{Solutions of the generalized eigenvalue problems in the differential realization}
We shall finally examine how solving the GEVP $\mathcal{M}f(z) = n \mathcal{L} f(z)$ and the adjoint problem compares to the representation theoretic computations that were performed of the overlaps $\widetilde{\langle z|}P_n(\alpha, \beta)\rangle$ and $\langle z |L^T |Q_m(\alpha, \beta)\rangle$. A first look shows that the GEVPs in the differential model will be of hypergeometric nature as is confirmed by the expressions of the overlaps. Let us focus on this more closely.

Given the expressions \eqref{L1} and \eqref{L2} for $\mathcal{L}$ and $\mathcal{M}$, we see that  $\mathcal{M}f(z) = n \mathcal{L} f(z)$ takes the form of the hypergeometric equation \cite{bateman1953higher}
\begin{equation}
    z (1-z) \frac{d^2 f}{dz^2} + \big[c - (a + b + 1)z\big] \frac{df}{dz} - ab f = 0, \label{hyper}
\end{equation}
with parameters
\begin{equation}
    a = -n, \qquad b = \alpha + 1, \qquad c= 1 - n - \beta. \label{par1}
\end{equation}
In the following we shall use Bateman's nomenclature \cite{bateman1953higher} of the $24$ Kummer solutions; these are arranged in six sets such that the four elements in each set represent the same function. The representatives $u_1, u_2, ..., u_6$ of the sets are in general different although \eqref{Kummer1} is a case where $u_1 \propto u_2$. With the parameters given by \eqref{par1}, it is immediate to see that the solution 
\begin{equation}
    u_1 = {_2}F_1 \left( {a, \; b\atop c} ; z\right)
\end{equation}
will yield directly (up to a constant) the Askey polynomials $P_n(z; \alpha, \beta)$. 

Consider now the adjoint operators:
\begin{equation}
\mathcal{L}^T = (1-z) \partial_z + \alpha \mathcal{I} \label{L1C} 
\end{equation}
\begin{equation}
\mathcal{M}^T = z(z-1) \partial_z^2 + \big[(2-\alpha)z -\beta-1 \big] \partial_z - \alpha \mathcal{I}. \label{L2C} 
\end{equation}
The adjoint GEVP $\mathcal{M}^T f^*(z) = m \mathcal{L} f^*(z)$ also turns out to yield the hypergeometric equation \eqref{hyper} but this time with parameters
\begin{equation}
    a = m + 1, \qquad b = - \alpha, \qquad c = 1 + \beta + m. \label{par2} 
\end{equation}
Recall that $\mathcal{L}^T f^*(z)$ will provide a solution orthogonal to $f(z)$. Selecting $u_1$ for $f^*$ also will lead to functions trivially orthogonal over the unit circle. Consider instead 
\begin{equation}
    f^*(z) = u_3 = (1 - z)^{-a} {_2}F_1 \left( {a, \;\;  c - b\atop a + 1 - b} \; ; \frac{1}{1 - z}\right).
\end{equation}
Using
\begin{equation}
    \frac{(k + m + \alpha + 1)}{(m + \alpha + 2)_k} = \frac{(m + \alpha + 1)}{(m + \alpha + 1)_k}
\end{equation}
it is easy to find that in that case,
\begin{equation}
    \mathcal{L}^T f^*(z) = ( m + \alpha + 1)(1-z)^{-m - 1} {_2}F_1 \left( {m + 1, \;m + \alpha + \beta + 1\atop m + \alpha + 1} ; \frac{1}{1 - z}\right).
\end{equation}
We thus see that choosing the  solution $u_3$,  yields the result \eqref{matQ} obtained algebraically for $\langle z |L^T |P_m (\alpha, \beta) \rangle$. (Recall that the constant $a$ in this expression is here set equal to $1$.) The reader is reminded of eq. \eqref{overlap} where this function is seen to be composed of two parts: one, a power series in $z$ and the other, the function orthogonal to the Askey polynomial dressed with the hypergeometric
weight. 

Let us point out that within the differential realization, it is possible to pick a solution of $\mathcal{M}^T f^*(z) = m \mathcal{L} f^*(z)$ that will solely give the biorthogonal partner $Q_m(\frac{1}{z}, \alpha, \beta)$ multiplied by the weight. Indeed take $f^*(z)$ to be rather given by the solution
\begin{equation}
    u_4 = (1 - z)^{-b} {_2}F_1 \left( {b, \;\;  c - a\atop b + 1 - a} \; ; \frac{1}{1 - z}\right).
\end{equation}
Substituting the parameters \eqref{par2} we have in this instance
\begin{align}
    f^*(z) &= (1 - z)^{\alpha} {_2}F_1 \left( {-\alpha, \; \beta\atop -\alpha - m} \; ; \frac{1}{1 - z}\right) \nonumber\\
    &= \sum_{k=0}^{\infty} \frac{(-\alpha)_k (\beta)_k}{(-\alpha - m)_k} \; \frac{(1 - z) ^{-k + \alpha}}{k!}.
\end{align}
The action of $\mathcal{L}^T$ is again readily computed when the argument is a function of $(1 - z)$:
\begin{align}
    \mathcal{L}^Tf^*(z) &= \sum_{k=1}^{\infty} \frac{(-\alpha)_k (\beta)_k}{(-\alpha - m)_k} \; \frac{(1 - z) ^{-k + \alpha}}{(k-1)!} \nonumber \\
    &=\sum_{l=0}^{\infty} \frac{(-\alpha)_{l+1} (\beta)_{l+1}}{(-\alpha - m)_{l+1}} \; \frac{(1 - z) ^{-l + \alpha -1}}{l!} \nonumber \\
    & = \frac{\alpha \beta}{(\alpha +m)} (1 - z)^{\alpha -1} {_2}F_1 \left( {-\alpha + 1, \; \beta + 1 \atop 1 -\alpha - m} \; ; \frac{1}{1 - z}\right),
\end{align}
where we have used $(x)_{l+1} = x (x + 1)_l$. We thus observe that the action of $\mathcal{L}^T$ is to effect : $ \alpha \rightarrow \alpha - 1, \; \beta \rightarrow \beta + 1$; in view of the parameter identification \eqref{par2}, the action of $\mathcal{L}^T$ on $u_4$ yields again $u_4$ with the following parameters:
\begin{equation}
    a = m + 1, \qquad b = -\alpha + 1, \qquad c = 2 + \beta + m. \label{par3}
\end{equation}
Now another expression for $u_4$ is 
\begin{equation}
    u_4 = (-z)^{a - c} (1 - z)^{c-a-b} {_2}F_1 \left( {1 - a, \;\;  c - a\atop b + 1 - a} \; ; \frac{1}{z}\right).
\end{equation}
Using the parameters \eqref{par3}, we then find that choosing $u_4$ as solution of the hypergeometric equation stemming from the adjoint GEVP $\mathcal{M}^T f^*(z) = m \mathcal{L} f^*(z)$ leads to 
\begin{equation}
    \mathcal{L}^Tf^*(z) \propto (-z) ^{-1 - \beta} (1 - z)^{\alpha + \beta} Q_m(\frac{1}{z}, \alpha, \beta).
\end{equation}
That is, we obtain as unique term, up to a factor, the orthogonal partner of $P_n(z; \alpha, \beta)$ multiplied by the weight.
\section{Conclusion}\label{C}
It is now time to wrap up and offer perspectives. We have presented a unified algebraic interpretation of the biorthogonal Askey polynomials on the circle and of the Jacobi polynomials on the interval $[0, 1]$. It is based on an algebra with three generators $L, M, X$ verifying quadratic relations that we have called the meta-Jacobi algebra and denoted by $m\mathfrak{J}$. The Askey polynomials $P_n(z; \alpha, \beta)$ arise as overlaps between the basis elements that are on the one hand the solutions, on an infinite dimensional module, of the generalized eigenvalue problem defined by the generators $L$ and $M$ and on the other hand, the eigenvectors of the adjoint of $X$. The biorthogonal partners $Q_n(z; \alpha, \beta)$ are obtained similarly from the reciprocal adjoints. The same framework was seen to provide an algebraic picture for the Jacobi polynomials as overlaps between the eigenbases of $M$ and of $X^T$ (or of $M^T$ and $X$). Proofs of the orthogonality relations were found to follow. With the introduction of a differential model for the meta-Jacobi algebra, the bispectrality of the Askey polynomials $P_n(z; \alpha, \beta)$ was accounted for in particular; their differential equation and the recurrence relation were explicitly obtained and found to be of GEVP form. 

The meta-Jacobi algebra is actually isomorphic to the Lie algebra $\mathfrak{su}(1, 1)$. We nevertheless kept with the (possibly redundant) terminology because the relevant presentation is parallel to that of the meta-Hahn algebra previously introduced \cite{tsujimoto2021algebraic}, \cite{vinet2020unified} to treat in a unified way orthogonal polynomials and biorthogonal rational functions of the Hahn type. Both algebras involve a non-commutative generalization of the plane which is supplemented by the addition of a third generator that brings to the fore significant representation theoretic features. The use of GEVPs proved to be a key aspect. As said already  in the Introduction, we see a pattern develop and we suspect that it might be possible to associate meta-algebras to most entries of the Askey scheme so as to simultaneously describe the bispectrality of the hypergeometric polynomials and of associated biorthogonal (rational) functions. This most likely relates to the forays by Kim and Stanton \cite{kim2020orthogonal} towards the development of a scheme for orthogonal polynomials of type $R_I$. To be sure, it is with enthusiasm that we plan to pursue the investigations of meta-algebras and their relations to special functions.

\appendix

\section{The computation of $\langle z | L^T | Q_m\rangle$} \label{details1}
 We provide in this Appendix the details on how the result of Proposition \ref{partner} is obtained.
 
 Given the expression \eqref{matQ} for $\langle z | L^T | Q_m\rangle$, we use the following linear relation between Kummer solutions of the hypergeometric equation \cite{bateman1953higher}(Sect 2.9, eq. (35)):
\begin{align}
    &{_2}F_1\left({a,\; b \atop c }; z \right) = \label{comeon}
    \frac{\Gamma (a + 1 -c) \Gamma (b + 1 - c)}{\Gamma (a + b + 1 - c) \Gamma (1 - c)}  \:{_2}F_1\left({a,\; \; b \atop a + b + 1 - c }; 1- z \right)  \\ \nonumber \\[0.5pt] \nonumber 
 & - \; \frac{\Gamma (a + 1 - c) \Gamma (b + 1 - c) \Gamma ( c - 1)}{\Gamma (a) \Gamma (b) \Gamma (1 - c)} \; z^{1-c} (1-z)^{c - a - b} \;   {_2}F_1\left({1 - a,\; 1 - b \atop 2 - c }; z \right).
\end{align}
This yields 
\begin{align}
    &{_2}F_1\left({m+1,\; m + \alpha + \beta + 1 \atop m + \alpha + 1 }; \frac{1}{1-z} \right) =  \nonumber \\[0.5pt]  \nonumber \\
    &\frac{\Gamma (1 - \alpha) \Gamma (\beta + 1)}{\Gamma (m + \beta + 2) \Gamma (- m - \alpha)}  \:{_2}F_1\left({m + 1,\; m + \alpha + \beta + 1 \atop m + \beta + 2 }; \frac{z}{z-1} \right) \; - \label{2f1}   \\[0.5pt] \nonumber \\
 & \frac{\Gamma (1 - \alpha) \Gamma (\beta + 1) \Gamma (m + \alpha)}{m! \; \Gamma (m + \alpha + \beta + 1) \Gamma (-m - \alpha)} \; (-z)^{-m - \beta - 1} (1-z)^{2m + \alpha + \beta + 1} \;  
 {_2}F_1\left({-m,\; - m - \alpha - \beta \atop 1 - m - \alpha }; \frac{1}{1-z} \right).\nonumber
\end{align}
Now use \cite{bateman1953higher} (Sect. 2.9 eqs. (1) \& (4) and (9) \& (11)):
\begin{equation}
    {_2}F_1\left({a,\; c - b \atop c }; \frac{z}{z - 1} \right) = (1 - z)^{a} \; {_2}F_1\left({a,\; b \atop c }; z \right)
\end{equation}
and 
\begin{equation}
    {_2}F_1\left({a,\; c - b \atop a + 1 - b }; \frac{1}{1-z} \right) = (1 - z)^{a} (-z)^{- a} \; {_2}F_1\left({a,\; a + 1 - c \atop a + 1 - b }; \frac{1}{z} \right)
\end{equation}
to reexpress the two ${_2}F_1$'s on the right hand side of \eqref{2f1} as functions of $z$ and $\frac{1}{z}$ respectively. Recalling then the definition \eqref{part} of the biorthogonal partner $Q_m(z, \alpha, \beta)$ of the Askey polynomials, one arrives at \eqref{overlap} with the help of the relation
\begin{equation}
    \Gamma (-m - \alpha) \Gamma (m + \alpha + 1)= (-1)^{m+1} \Gamma (\alpha) \Gamma (1 - \alpha) \label{useful}
\end{equation}
which is a consequence of the identity
\begin{equation}
    \Gamma(x) \Gamma (1-x) = \frac{\pi}{\sin \;\pi x}. \label{gammaid}
\end{equation}

\section{The determination of $\langle z \widetilde{| J_n \rangle}$}\label{compJac} 

Details on how formula \eqref{overlap2} for $\langle z \widetilde{| J_m \rangle}$ is obtained are given here. We need to transform the ${_2}F_1$ that occurs in the expression \eqref{secondoverlap} of this overlap. First we use the following relation between three solutions of the hypergeometric equation \cite{bateman1953higher} (Sect 2.9, eq. (34)):
\begin{align}
    &{_2}F_1\left({a,\; b \atop c }; z \right) = \label{comeonagain}
    \frac{\Gamma (c) \Gamma (b - a)}{\Gamma (c - a) \Gamma (b)} (-z)^{-a} \:{_2}F_1\left({a,\; \; a + 1 - c \atop a + 1 - b }; \frac{1}{z} \right)  \\ \nonumber \\[0.5pt] \nonumber 
 & + \;\frac{\Gamma (c) \Gamma (a - b)}{\Gamma (c-b) \Gamma (a)} \; (-z)^{a-c} (1-z)^{c - a - b} \;   {_2}F_1\left({1 - a,\; c - a \atop b + 1 - a }; \frac{1}{z} \right).
\end{align}
From this identity we find:
\begin{align}
    &{_2}F_1\left({m+1,\; m + \alpha + \beta + 1 \atop 2m + \alpha + 2 }; \frac{1}{1-z} \right) =  (z - 1 )^{m+1} \nonumber \\[0.5pt]  \nonumber \\
    &\times \;\Bigg[\frac{\Gamma (2m + \alpha + 2) \Gamma (\alpha + \beta)}{\Gamma (m + \alpha + 1) \Gamma (m + \alpha + \beta + 1)}  \:{_2}F_1\left({m + 1,\; -m - \alpha \atop 1 - \alpha - \beta }; 1 - z \right) \; + \label{cheers}   \\[0.5pt] \nonumber \\
 & \frac{\Gamma (2m + \alpha + 2) \Gamma (- \alpha - \beta)}{m! \;\Gamma (m + 1 - \beta)} \; z^{-\beta}(z - 1)^{\alpha + \beta} \;  
 {_2}F_1\left({-m,\; m + \alpha + 1 \atop \alpha + \beta + 1}; 1-z \right)\Bigg].\nonumber 
\end{align}
We now apply the relation \eqref{comeon} to convert each of the two ${_2}F_1$s on the right hand side of \eqref{cheers} that are functions of $(1-z)$ into 
combinations of ${_2}F_1$s that are functions of $z$. This leads to
\begin{align}
    &{_2}F_1\left({m+1,\; m + \alpha + \beta + 1 \atop 2m + \alpha + 2 }; \frac{1}{1-z} \right) = \nonumber \\[0.5pt]  \nonumber \\ 
    &\frac{\Gamma (2m + \alpha + 2) \Gamma(\alpha + \beta)\Gamma (1 - \alpha - \beta) \Gamma(-\beta) }{\Gamma (m + \alpha + 1) \Gamma (m + \alpha + \beta + 1) \Gamma (1 +m - \beta) \Gamma (-m - \alpha - \beta) } (z-1)^{m + 1} \:{_2}F_1\left({m + 1,\; -m - \alpha \atop 1 + \beta }; z \right) \nonumber \\ \nonumber\\[5pt] 
   & + \Bigg[(-1)^{(\alpha + \beta)}\frac{\Gamma (2m + \alpha + 2) \Gamma (\alpha +\beta)\Gamma (1 - \alpha - \beta) \Gamma (\beta)}{m! \; \Gamma (m + \alpha  + 1) \Gamma (-m - \alpha) \Gamma (m + \alpha + \beta +1)}\label{cheerios}   \\[0.5pt] \nonumber \\
   & + \frac{\Gamma( 2m + \alpha + 2) \Gamma (- \alpha - \beta) \Gamma (\alpha + \beta + 1) \Gamma (\beta)}{m! \;\Gamma (m + 1 - \beta)\Gamma (-m + \beta) \Gamma (m + \alpha + \beta + 1)}\Bigg]
  z^{-\beta}(z - 1)^{m + 1 +\alpha + \beta} \;  
 {_2}F_1\left({-m,\; m + \alpha + 1 \atop 1 - \beta}; z \right). \nonumber
\end{align}
Simplifications are carried out through repeated use of the identity \eqref{gammaid} and various implications such as \eqref{useful} and by observing that
\begin{align}
    &\frac{1}{\sin\; \pi(\alpha + \beta)} \left((-1)^{\alpha + \beta} \sin\; \pi \alpha + \sin\; \pi \beta  \right) \nonumber \\
    &= \frac{e^{i \pi \alpha}}{\sin\; \pi (\alpha + \beta)} \left(e^{i \pi \beta} \sin\; \pi \alpha + e^{-i \pi \alpha} \sin\; \pi \beta\right) 
    = e^{i \pi \alpha} = (-1)^{\alpha}. \label{littleid}
\end{align}
One finally obtains
\begin{align}
    &{_2}F_1\left({m+1,\; m + \alpha + \beta + 1 \atop 2m + \alpha + 2 }; \frac{1}{1-z} \right) =  (1 - z )^{m+1} \nonumber \\[0.5pt]  \nonumber \\
    &\times \;\Bigg[\frac{\Gamma (2m + \alpha + 2) \Gamma(-\beta)}{\Gamma (m + \alpha + 1) \Gamma (m - \beta + 1)}  \:{_2}F_1\left({m + 1,\; -m - \alpha \atop 1 + \beta }; z \right) \; + \\[0.5pt] 
    &\frac{\Gamma (2m + \alpha + 2) \Gamma (\beta)}{m! \;\Gamma (m + \alpha + \beta +1)}\label{cheerios} \nonumber(-z)^{-\beta}(1-z)^{\alpha + \beta} \;  
 {_2}F_1\left({-m,\; m + \alpha + 1 \atop 1 - \beta}; z \right)\Bigg]
\end{align}
which readily gives \eqref{overlap2}.
\section{Negative indices}\label{neg}

In the main part of the paper, it sufficed for the purpose of interpreting the Askey polynomials and their biorthogonal partners to focus on GEVP and EVP solutions with non-negative (integer) eigenvalues. For completeness we briefly indicate in this appendix how the situations with negative integers can be treated and seen to lead to redundant information.

\subsection{}
Consider equation \eqref{recurd} and assume that $n < 0$. Let 
\begin{equation}
    n = - s - 1 , \qquad s = 0, 1, \dots \label{n}
\end{equation}
In this case the recursion relation still implies $d_n(k) = 0$ for $k > n$ but no longer bounds $k$ from below. Write $k$ in the form
\begin{equation}
    k = - s - 1 - l, \qquad l = 0, 1, \dots \label{k}
\end{equation}
Upon substituting \eqref{n} and \eqref{k} and taking $d_n(k) \equiv \Tilde{d}_s(s+l)$, equation \eqref{recurd} becomes
\begin{equation}
    (s + l) (l + s - \alpha - \beta) \Tilde{d}_s (s + l - 1) + l (l + s - \alpha) \Tilde{d}_s (s + l) = 0.
\end{equation}
We observe that this last relation coincide with the condition \eqref{recur*} that had been obtained from the adjoint GEVP with a positive eigenvalue under the substitutions
\begin{equation}
    m \rightarrow s, \quad \alpha \rightarrow - \alpha - 1, \quad \beta \rightarrow - \beta + 1, \quad d_m ^*(m+l) \rightarrow  \Tilde{d}_s(s+l).
\end{equation}
Hence,
\begin{equation}
    d_n(k)= d_{-s-1}(-s-1-l)=(-1)^l \frac{(s+1)_l (s-\alpha - \beta + 1)_l}{l! (s - \alpha + 1)_l} \; d_n(n), \qquad s, l = 0, 1, 2, \dots.
\end{equation}
\subsection{}
Examine now equation \eqref{recur*} when $m < 0$. Set $m=-s-1, s=0, 1, \dots$. In this case the recursion equation implies that $d^*_m (k)=0$ for $k<m$ and also truncates at $k=0$. The non-zero values of $d_n^*(k)$ therefore only occur for
\begin{equation}
    k = - l - 1, \qquad l=0, \dots s.
\end{equation}
Incorporating the above redefinitions in \eqref{recur*} and taking $d^*_m(k)=d_{-s-1}^*(-l-1) \equiv \Tilde{d}_m^*(l)$ we get
\begin{equation}
    (l+1)(l-\alpha - \beta + 1) \Tilde{d}_m^*(l+1) + (l-s)(l-\alpha) \Tilde{d}_m^*(l) = 0
\end{equation}
and see that this equation can be retrieved from \eqref{recurd} under the substitutions
\begin{equation}
    n \rightarrow m, \quad \alpha \rightarrow -\alpha - 1, \quad \beta \rightarrow - \beta - 1, \quad d_n(k) \rightarrow \Tilde{d}_m^*(l).
\end{equation}
It follows that for negative $m$
\begin{equation}
    d_m^*(k)=d_{-s-1}^*(-l - 1) = (-1)^l \frac{(-s)_l(-\alpha)_l}{l! (-\alpha -\beta + 1)_l} d^*_n(-1), \qquad l=0, \dots s.
\end{equation}
\subsection{}
We may check the orthogonality of $|P_n\rangle$ and $L^T|Q_m\rangle$, $m \neq n$, for various possibilities regarding the sign of the indices $m$ and $n$. In summary, the summation ranges are as follows:
\begin{itemize}
    \item For $n \geq 0, m\geq 0$
    \begin{align}
        &|P_n\rangle = \sum_{k=0} ^{n} d_n(k) |\tau, k\rangle, \\
        &|Q_m\rangle = \sum_{k=m} ^{\infty} d_m(k) |\tau, k\rangle;
    \end{align}
    \item For $n < 0, m < 0$
        \begin{align}
        &|P_n\rangle = \sum_{k=-\infty} ^{n} d_n(k) |\tau, k\rangle, \\
        &|Q_m\rangle = \sum_{k=m} ^{-1} d_m(k) |\tau, k\rangle.
    \end{align}
\end{itemize}
It is manifest that the orthogonality prevails when one index is non-negative and the other is negative. When the two indices are negative, the proof of orthogonality follows the one given for two non-negative indices since as we observed the change of signs basically flips the coefficients $d$ and $d^*$.

\subsection{}
Regarding the special functions, in light of this exchange of the expansion coefficients, the roles of $|P_n\rangle$ and of $L^T|Q_m\rangle$ are inverted when the indices are negative. For instance, we have
\begin{equation}
    |Q_{-s-1}\rangle = \sum_{l=0}^s (-1)^l \frac{(-s)_l(-\alpha)_l}{l! (-\alpha - \beta + 1)_l} |\tau, -l-1\rangle.
\end{equation}
The overlap of $L^T|Q_{-s-1}\rangle$ with the state $|z\rangle$ given in
\eqref{ez} is then found to be
\begin{equation}
    \langle z |L^T|Q_{-s-1}\rangle = \alpha d_{-s-1}(-1) (z-1)^{\Tilde{a}-1} {_2}F_1\left({-s,\; 1-\alpha \atop 1-\alpha -\beta }; 1- z \right).
\end{equation}
Owing again to \eqref{Kummer1}, we see that the Askey polynomials arise in this case in the overlap $ \langle z |L^T|Q_{-s-1}\rangle$ with a change of parameters.
\subsection{}
Things can be seen to proceed similarly in the treatment of the Jacobi polynomials if negative indices are considered.
\section*{Acknowledgments}
The authors are  grateful to Tom Koornwinder for correspondence and bringing some references to their attention. They have much appreciated Erik Koelink's comments on the manuscript and are  thankful to Julien Gaboriaud for kind assistance. The work of LV is supported in part by a Discovery Grant from the Natural Sciences and Engineering Research Council (NSERC) of Canada. AZ who is funded by the National Foundation of China (Grant No.11771015) gratefully acknowledges the hospitality of the CRM over an extended period and the award of a Simons CRM professorship..
\bibliographystyle{unsrt} 
\bibliography{ref_AJ.bib}

\end{document}